\documentclass{amsart}
\usepackage{amssymb}
\usepackage{amsmath}

\usepackage{amsfonts}

\newtheorem{theorem}{Theorem}
\theoremstyle{plain}

\newtheorem{definition}{Definition}

\newtheorem{lemma}{Lemma}

\newtheorem{proposition}{Proposition}
\newtheorem{remark}{Remark}

\numberwithin{equation}{section}
\input{tcilatex}

\begin{document}
\title[]{Homogenization of a linear elastic body with rigid inclusions and a
Robin type boundary conditions}
\author{Lazarus Signing}
\address{University of Ngaoundere, Department of Mathematics and Computer
Science, P.O.Box 454 Ngaoundere (Cameroon)}
\email{lsigning@yahoo.fr}
\dedicatory{}
\thanks{}

\begin{abstract}
This paper is devoted to study of the limiting behaviour of an elastic
material with periodically distributed rigid inclusions of size $\varepsilon 
$, as the small parameter $\varepsilon $ goes to zero. We address here the
case with inclusions of the same size as the period of the structure. The
body in consideration here is suppose to be clamped on one part of its
exterior boundary and submitted to given tractions on the other. By means of
the well known two-scale convergence techniques, one convergence result is
proved.
\end{abstract}

\maketitle

\section{Introduction}

Let $\Omega $ be a smooth bounded open set of $\mathbb{R}_{x}^{N}$ (the $N$%
-dimensional numerical space of variables $x=\left( x_{1},...,x_{N}\right) $%
) with $N\geq 2$. We consider a compact subset $T$ of $\mathbb{R}_{y}^{N}$
with smooth boundary and nonempty interior such that 
\begin{equation}
T\subset Y=\left( -\frac{1}{2},\frac{1}{2}\right) ^{N}\text{.}  \label{eq1.1}
\end{equation}%
For any $\varepsilon >0$, we define 
\begin{equation}
t^{\varepsilon }=\left\{ k\in \mathbb{Z}^{N}:\varepsilon \left( k+T\right)
\subset \Omega \right\} \text{,}  \label{eq1.1a}
\end{equation}%
\begin{equation}
T^{\varepsilon }=\underset{k\in t^{\varepsilon }}{\cup }\varepsilon \left(
k+T\right)  \label{eq1.1b}
\end{equation}%
and%
\begin{equation}
\Omega ^{\varepsilon }=\Omega \backslash T^{\varepsilon }\text{,}
\label{eq1.1c}
\end{equation}%
where $\mathbb{Z}$ denotes the integers. Throughout this study, $\Omega
^{\varepsilon }$ is a medium constituted of an elastic porous material, and $%
T$ is the reference rigid part while $\varepsilon \left( k+T\right) $ is a
rigid obstacle of size $\varepsilon $. The subset $T^{\varepsilon }$ is
therefore the union of rigid particules of size $\varepsilon $ in the porous
domain $\Omega ^{\varepsilon }$. We denote by $\mathbf{n}=\left(
n_{j}\right) _{1\leq j\leq N}$ the outward unit normal to $\partial
T^{\varepsilon }$ with respect to $\Omega ^{\varepsilon }$.

For any Roman character such as $i$, $j$ (with $1\leq i,j\leq N$), $u^{i}$
(resp. $u^{j}$) denotes the $i$-th (resp. $j$-th) component of a vector
function $\mathbf{u}$ in $L_{loc}^{1}\left( \Omega \right) ^{N}$ or in $%
L_{loc}^{1}\left( \mathbb{R}^{N}\right) ^{N}$. Further, for any real $%
0<\varepsilon <1$, we define $u^{\varepsilon }$ as 
\begin{equation*}
u^{\varepsilon }\left( x\right) =u\left( \frac{x}{\varepsilon }\right) \text{%
\qquad }\left( x\in \Omega \right)
\end{equation*}%
for $u\in L_{loc}^{1}\left( \mathbb{R}_{y}^{N}\right) $. More generally, for 
$u\in L_{loc}^{1}\left( \Omega \times \mathbb{R}_{y}^{N}\right) $, it is
customary to put%
\begin{equation*}
u^{\varepsilon }\left( x\right) =u\left( x,\frac{x}{\varepsilon }\right) 
\text{\qquad }\left( x\in \Omega \right)
\end{equation*}%
whenever the right-hand side makes sense (see, e.g., \cite{bib8}).

Let $a_{ijkh}$ $\left( 1\leq i,j,k,h\leq N\right) $ and $\mathcal{\theta }$
be real functions in $L^{\infty }\left( \mathbb{R}^{N}\right) $ such that: 
\begin{equation}
a_{ijkh}=a_{jihk}=a_{khij}\text{, and }\mathcal{\theta }\left( y\right) \geq
\alpha _{0}\text{ a.e.in }y\in \mathbb{R}^{N}\text{,}  \label{eq1.2a}
\end{equation}%
\begin{equation}
\sum_{i,j;k;h=1}^{N}a_{ijkh}\left( y\right) \zeta _{ij}\zeta _{kh}\geq
\alpha \left\vert \zeta \right\vert ^{2}\text{\qquad }\left( \zeta =\left(
\zeta _{ij}\right) \in \mathbb{R}^{N\times N}\right) \text{ a.e. in }y\in 
\mathbb{R}^{N}\text{,}  \label{eq1.2}
\end{equation}%
where $\alpha >0$ and $\alpha _{0}>0$ are constants. Let us denote by $%
\mathbf{e}$ and $\mathbf{\sigma }^{\varepsilon }$ the strain and stress
tensor respectively, related in the framework of linear elasticity by the
Hooke law:%
\begin{equation*}
\mathbf{\sigma }_{ij}^{\varepsilon }=\sum_{k,h=1}^{N}a_{ijkh}^{\varepsilon }%
\mathbf{e}_{kh}\left( \mathbf{u}_{\varepsilon }\right)
\end{equation*}%
where 
\begin{equation*}
\mathbf{e}_{ij}\left( \mathbf{u}_{\varepsilon }\right) =\frac{1}{2}\left( 
\frac{\partial u_{\varepsilon }^{i}}{\partial x_{j}}+\frac{\partial
u_{\varepsilon }^{j}}{\partial x_{i}}\right) \qquad \left( 1\leq i,j\leq
N\right)
\end{equation*}%
for $\mathbf{u}_{\varepsilon }=\left( u_{\varepsilon }^{j}\right) $. Let us
also suppose that $\Gamma =\overline{\Gamma }_{1}\cup \overline{\Gamma }_{2}$%
, $\Gamma _{1}$ and $\Gamma _{2}$ being the disjoint open parts of the
smooth boundary $\Gamma $ of $\Omega $ such that $meas\left( \Gamma
_{1}\right) >0$.

For any fixed $0<\varepsilon <1$, we consider the boundary value problem%
\begin{equation}
\quad \quad -\func{div}\mathbf{\sigma }^{\varepsilon }=\mathbf{f}%
^{\varepsilon }\text{ in }\Omega ^{\varepsilon }\text{,}  \label{eq1.3}
\end{equation}%
\begin{equation}
\quad \qquad \qquad \mathbf{u}_{\varepsilon }=0\text{ on }\Gamma _{1}\text{,}
\label{eq1.4}
\end{equation}%
\begin{equation}
\qquad \qquad \mathbf{\sigma }^{\varepsilon }\mathbf{n}=\mathbf{t}\text{ on }%
\Gamma _{2}\text{,}  \label{eq1.5}
\end{equation}%
\begin{equation}
\text{\qquad \qquad \qquad \quad }\mathbf{\sigma }^{\varepsilon }\mathbf{n}%
=-\varepsilon \mathcal{\theta }^{\varepsilon }\mathbf{u}_{\varepsilon }\text{
on }\partial T^{\varepsilon }  \label{eq1.7}
\end{equation}%
where $\mathbf{f}=\left( f_{j}\right) \in L^{\infty }\left( \mathbb{R}%
_{y}^{N}\right) ^{N}$ and $\mathbf{t=}\left( t_{j}\right) \in L^{2}\left(
\Gamma \right) ^{N}$ are vector functions with real components. For the
variational formulation of (\ref{eq1.3})-(\ref{eq1.7}), let us introduce%
\begin{equation*}
\mathbf{V}_{\varepsilon }=\left\{ \mathbf{v}\in H^{1}\left( \Omega
^{\varepsilon };\mathbb{R}\right) ^{N}:\mathbf{v=}0\text{ on }\Gamma
_{1}\right\} \text{,}
\end{equation*}%
where $H^{1}\left( \Omega ^{\varepsilon };\mathbb{R}\right) $ is the space
of\ functions in\ the Sobolev$\ $space$\ H^{1}\left( \Omega ^{\varepsilon
}\right) $ with real values, and let $\mathbf{a}^{\varepsilon }\left(
.,.\right) $ be the bilinear form on $H^{1}\left( \Omega ^{\varepsilon };%
\mathbb{R}\right) ^{N}$ given\ by%
\begin{equation*}
\mathbf{a}^{\varepsilon }\left( \mathbf{u},\mathbf{v}\right)
=\sum_{i,j,k,h=1}^{N}\int_{\Omega ^{\varepsilon }}a_{ijkh}^{\varepsilon }%
\mathbf{e}_{ij}\left( \mathbf{u}\right) \mathbf{e}_{kh}\left( \mathbf{v}%
\right) dx+\varepsilon \int_{\partial T^{\varepsilon }}\mathcal{\theta }%
^{\varepsilon }\mathbf{u}{\small \cdot }\mathbf{v}d\sigma _{\varepsilon }
\end{equation*}%
for $\mathbf{u=}\left( u^{k}\right) $ and $\mathbf{v=}\left( v^{k}\right)
\in H^{1}\left( \Omega ^{\varepsilon };\mathbb{R}\right) ^{N}$, the dot
denoting the Euclidean inner product, and $d\sigma _{\varepsilon }$ being
the surface measure on $\partial T^{\varepsilon }$. The boundary value
problem\ (\ref{eq1.3})-(\ref{eq1.7}) naturally\ implies the following
variational equation:%
\begin{equation}
\left\{ 
\begin{array}{c}
\mathbf{u}_{\varepsilon }\in \mathbf{V}_{\varepsilon },\qquad \qquad \qquad
\qquad \qquad \qquad \qquad \qquad \qquad \\ 
\mathbf{a}^{\varepsilon }\left( \mathbf{u}_{\varepsilon },\mathbf{v}\right)
=\int_{\Omega ^{\varepsilon }}\mathbf{f}^{\varepsilon }{\small \cdot }%
\mathbf{v}dx+\int_{\Gamma _{2}}\mathbf{t}{\small \cdot }\mathbf{v}d\Gamma 
\text{ for all }\mathbf{v}\in \mathbf{V}_{\varepsilon }\text{,}%
\end{array}%
\right.  \label{eq1.8}
\end{equation}%
$d\Gamma $ being the surface measure on $\Gamma $. The variational\ problem\
(\ref{eq1.8}) is a classical\ one which\ admits a unique solution, in view
of\ (\ref{eq1.2a})-(\ref{eq1.2}). Further, it is easy\ to check that (\ref%
{eq1.8}) leads to\ (\ref{eq1.3})-(\ref{eq1.7}). Thus, the problem (\ref%
{eq1.3})-(\ref{eq1.7}) admits a unique solution $\mathbf{u}_{\varepsilon }$\
in $\mathbf{V}_{\varepsilon }$.

Our aim here is to investigate the asymptotic behaviour, as $\varepsilon
\rightarrow 0$, of $\mathbf{u}_{\varepsilon }$ and $\mathbf{\sigma }%
^{\varepsilon }$ under the hypotheses that 
\begin{equation}
a_{ijkh}\left( y+k\right) =a_{ijkh}\left( y\right) \text{,\quad }\mathcal{%
\theta }\left( y+k\right) =\mathcal{\theta }\left( y\right) \text{ and }%
\mathbf{f}\left( y+k\right) =\mathbf{f}\left( y\right) \quad \left( 1\leq
i,j,k,h\leq N\right)  \label{eq1.8a}
\end{equation}%
for almost all $y\in \mathbb{R}^{N}$ and for all $k\in \mathbb{Z}^{N}$.

The study of this problem turns out to be of benefit to the modelling of an
heterogeneous elastic material with rigid periodically distributed
inclusions.

Many authors have addressed similar problems in several contexts, using
various methods. The homogenization of an elastic material with inclusions
in frictionless contact has been studied by MiKeli\'{c}, Shillor and Tapi%
\`{e}ro in \cite{bib6}. A viscoelastic periodically perforated material with
rigid inclusions with contact and friction described by linear conditions,
has been considered in \cite{bib4} by Gilbert, Panchenko and Xie. Further,
we mention the paper by Iosif'yan \cite{bib5} in which a system of linear
elasticity has been considered for a periodically perforated domain with a
nonlinear Robin condition on the boundary of the inclusions. In \cite{bib3},
by the \textit{periodic unfolding }method\textit{,} Capatina and Timofte
have addressed a similar problem with several linear and nonlinear
conditions on the boundary of the inclusions.

In this work, we consider a linear case of the condition studied in the
first problem of \cite{bib3}. In this case, the body is clamped on the face $%
\Gamma _{1}$, and submitted to external volume forces of density $\mathbf{f}%
^{\varepsilon }$ and a surface traction of density $\mathbf{t}$ on the rest
of the boundary $\Gamma _{2}$. On the boundary of the rigid inclusions $%
\partial T^{\varepsilon }$, we have a linear Robin type condition.

Our approach is the \textit{two-scale convergence }method, which is nothing
but the \textit{sigma-convergence }in the periodic setting. The results on
the two-scale convergence for periodic surfaces\textit{\ }in \cite{bib2}
(see also \cite{bib12}) are the keystone in this work and make easier a
rigorous proof of the convergence of the homogenization process. By this
means, we derive the macroscopic homogenized model for (\ref{eq1.3})-(\ref%
{eq1.7}).

Unless otherwise specified, vector spaces throughout are considered over the
complex field, $\mathbb{C}$, and scalar functions are assumed to take
complex values. Let us recall some basic notations. If $X$ and $F$ denote a
locally compact space and a Banach space, respectively, then we write $%
\mathcal{C}\left( X;F\right) $ for continuous mappings of $X$ into $F$, and $%
\mathcal{B}\left( X;F\right) $ for those mappings in $\mathcal{C}\left(
X;F\right) $ that are bounded. We denote by $\mathcal{K}\left( X;F\right) $
the mappings in $\mathcal{C}\left( X;F\right) $ having compact supports.\ We
shall assume $\mathcal{B}\left( X;F\right) $ to be equipped with the
supremum norm $\left\Vert u\right\Vert _{\infty }=\sup_{x\in X}\left\Vert
u\left( x\right) \right\Vert $ ($\left\Vert {\small \cdot }\right\Vert $
denotes the norm in $F$). For shortness we will write $\mathcal{C}\left(
X\right) =\mathcal{C}\left( X;\mathbb{C}\right) $, $\mathcal{B}\left(
X\right) =\mathcal{B}\left( X;\mathbb{C}\right) $ and $\mathcal{K}\left(
X\right) =\mathcal{K}\left( X;\mathbb{C}\right) $. Likewise in the case when 
$F=\mathbb{C}$, the usual spaces $L^{p}\left( X;F\right) $ and $%
L_{loc}^{p}\left( X;F\right) $ ($X$ provided with a positive Radon measure)
will be denoted by $L^{p}\left( X\right) $ and $L_{loc}^{p}\left( X\right) $%
, respectively. Finally, the numerical space $\mathbb{R}^{N}$ and its open
sets are each provided with Lebesgue measure denoted by $dx=dx_{1}...dx_{N}$.

The rest of the paper is organized as follows. Section 2 is devoted to the
preliminaries while in Section 3, a convergence result is proved for (\ref%
{eq1.3})-(\ref{eq1.7}).

\section{Preliminary results}

Before we begin with preliminaries, let us note that, if $\mathbf{w=}\left(
w^{k}\right) _{1\leq k\leq N}$ with $w^{k}\in L^{p}\left( \mathcal{O}\right) 
$, or if $\mathbf{w=}\left( w^{ij}\right) _{1\leq i,j\leq N}$ with $%
w^{ij}\in L^{p}\left( \mathcal{O}\right) $, where $\mathcal{O}$ is an open
set in $\mathbb{R}^{N}$, we will sometimes write $\left\Vert \mathbf{w}%
\right\Vert _{L^{p}\left( \mathcal{O}\right) }$ for $\left\Vert \mathbf{w}%
\right\Vert _{L^{p}\left( \mathcal{O}\right) ^{N}}$ or for $\left\Vert 
\mathbf{w}\right\Vert _{L^{p}\left( \mathcal{O}\right) ^{N\times N}}$.

Let us first recall the following result on the construction (for $%
\varepsilon >0$) of a suitable extension operator sending $H^{1}\left(
\Omega ^{\varepsilon };\mathbb{R}\right) ^{N}$ into $H^{1}\left( \Omega ;%
\mathbb{R}\right) ^{N}$.

\begin{proposition}
\label{pr2.1} For each real $\varepsilon >0$, there exists an operator $%
\mathcal{P}_{\varepsilon }$ of $H^{1}\left( \Omega ^{\varepsilon };\mathbb{R}%
\right) ^{N}$ into $H^{1}\left( \Omega ;\mathbb{R}\right) ^{N}$ with the
following properties: 
\begin{equation}
\mathcal{P}_{\varepsilon }\text{ sends continuously and linearly }%
H^{1}\left( \Omega ^{\varepsilon };\mathbb{R}\right) ^{N}\text{ into }%
H^{1}\left( \Omega ;\mathbb{R}\right) ^{N}\text{;}  \label{eq2.1}
\end{equation}%
\begin{equation}
\left( \mathcal{P}_{\varepsilon }\mathbf{v}\right) {\small \mid }_{\Omega
^{\varepsilon }}=\mathbf{v}\text{ for all }\mathbf{v}\in H^{1}\left( \Omega
^{\varepsilon };\mathbb{R}\right) ^{N}\text{;}  \label{eq2.2}
\end{equation}%
\begin{equation}
\left\Vert \mathcal{P}_{\varepsilon }\mathbf{v}\right\Vert _{L^{2}\left(
\Omega \right) }\leq c\left\Vert \mathbf{v}\right\Vert _{L^{2}\left( \Omega
^{\varepsilon }\right) }  \label{eq2.3a}
\end{equation}%
and%
\begin{equation}
\left\Vert \mathbf{e}\left( \mathcal{P}_{\varepsilon }\mathbf{v}\right)
\right\Vert _{L^{2}\left( \Omega \right) }\leq c\left\Vert \mathbf{e}\left( 
\mathbf{v}\right) \right\Vert _{L^{2}\left( \Omega ^{\varepsilon }\right) }
\label{eq2.3b}
\end{equation}
\ for all $\mathbf{v}\in H^{1}\left( \Omega ^{\varepsilon };\mathbb{R}%
\right) ^{N}$, where the constant $c>0$ depends solely on $Y$ and $T$.
\end{proposition}

The proof of the preceding proposition is to be found in, e.g., \cite[%
Theorem 4.2]{bib11}.

Let us go to our next purpose. We set 
\begin{equation}
\Theta =\tbigcup_{k\in \mathbb{Z}^{N}}\left( k+T\right) \text{.}
\label{eq2.6a}
\end{equation}%
Using the compactness of $T$, it is easy to check that $\Theta $ is closed
in $\mathbb{R}^{N}$. Let us set for any $\varepsilon >0$ 
\begin{equation}
Q^{\varepsilon }=\Omega \backslash \varepsilon \Theta \text{.}
\label{eq2.6b}
\end{equation}%
Then $Q^{\varepsilon }$ is an open set of $\mathbb{R}^{N}$ and clearly, $%
Q^{\varepsilon }\subset \Omega ^{\varepsilon }$. The set $Q^{\varepsilon }$
is made of two types of solid particules: on one hand, the solids of $\Omega
^{\varepsilon }$, on the other hand, the solids $\Omega \cap \varepsilon
\left( k+T\right) $ where $\varepsilon \left( k+T\right) $ intersects $%
\partial \Omega $.

Now, let us turn to some fundamental preliminary results on the \textit{%
sigma-convergence }in the periodic setting.

Let us first recall that a function $u\in L_{loc}^{1}\left( \mathbb{R}%
_{y}^{N}\right) $ is said to be $Y$-periodic if for each $k\in \mathbb{Z}%
^{N} $, we have $u\left( y+k\right) =u\left( y\right) $ almost everywhere
(a.e.) in $y\in \mathbb{R}^{N}$. If in addition $u$ is continuous, then the
preceding equality holds for every $y\in \mathbb{R}^{N}$, of course. The
space of all $Y$-periodic continuous complex functions on $\mathbb{R}%
_{y}^{N} $ is denoted by $\mathcal{C}_{per}\left( Y\right) $; that of all $Y$%
-periodic functions in $L_{loc}^{p}\left( \mathbb{R}_{y}^{N}\right) $ $%
\left( 1\leq p<\infty \right) $ is denoted by $L_{per}^{p}\left( Y\right) $. 
$\mathcal{C}_{per}\left( Y\right) $ is a Banach space under the supremum
norm on $\mathbb{R}^{N}$, whereas $L_{per}^{p}\left( Y\right) $ is a Banach
space under the norm 
\begin{equation*}
\left\Vert u\right\Vert _{L^{p}\left( Y\right) }=\left( \int_{Y}\left\vert
u\left( y\right) \right\vert ^{p}dy\right) ^{\frac{1}{p}}\text{ }\left( u\in
L_{per}^{p}\left( Y\right) \right) \text{.}
\end{equation*}

We will need the space $H_{per}^{1}\left( Y\right) $ of functions in $%
H_{loc}^{1}\left( \mathbb{R}_{y}^{N}\right) =W_{loc}^{1,2}\left( \mathbb{R}%
_{y}^{N}\right) $ which are $Y$-periodic, and the space $H_{\#}^{1}\left(
Y\right) $ of functions $u\in H_{per}^{1}\left( Y\right) $ such that $%
\int_{Y}\left( y\right) dy=0$. Provided with the gradient norm, 
\begin{equation*}
\left\Vert u\right\Vert _{H_{\#}^{1}\left( Y\right) }=\left(
\int_{Y}\left\vert \nabla _{y}u\right\vert ^{2}dy\right) ^{\frac{1}{2}}\text{
}\left( u\in H_{\#}^{1}\left( Y\right) \right) \text{,}
\end{equation*}%
where $\nabla _{y}u=\left( \frac{\partial u}{\partial y_{1}},...,\frac{%
\partial u}{\partial y_{N}}\right) $, $H_{\#}^{1}\left( Y\right) $ is a
Hilbert space.

Before we can recall the concept of \textit{sigma-convergence} in the
present periodic setting or the \textit{two-scale convergence}, let us
introduce one further notation. The letter $E$ throughout will denote a
family of real numbers $0<\varepsilon <1$ admitting $0$ as an accumulation
point. For example, $E$ may be the whole interval $\left( 0,1\right) $; $E$
may also be an ordinary sequence $\left( \varepsilon _{n}\right) _{n\in 
\mathbb{N}}$ with $0<\varepsilon _{n}<1$ and $\varepsilon _{n}\rightarrow 0$
as $n\rightarrow \infty $. In the latter case $E$ will be referred to as a 
\textit{fundamental sequence}.

Let $\Omega $ be a bounded open set in $\mathbb{R}_{x}^{N}$ and let $1\leq
p<\infty $.

\begin{definition}
\label{def2.1} A sequence $\left( u_{\varepsilon }\right) _{\varepsilon \in
E}\subset L^{p}\left( \Omega \right) $ is said to:

(i) weakly $\Sigma $-converge in $L^{p}\left( \Omega \right) $ to some $%
u_{0}\in L^{p}\left( \Omega ;L_{per}^{p}\left( Y\right) \right) $ if as

\noindent $E\ni \varepsilon \rightarrow 0$, 
\begin{equation}
\int_{\Omega }u_{\varepsilon }\left( x\right) \psi ^{\varepsilon }\left(
x\right) dx\rightarrow \int \int_{\Omega \times Y}u_{0}\left( x,y\right)
\psi \left( x,y\right) dxdy  \label{eq2.9}
\end{equation}%
\begin{equation*}
\begin{array}{c}
\text{for all }\psi \in L^{p^{\prime }}\left( \Omega ;\mathcal{C}%
_{per}\left( Y\right) \right) \text{ }\left( \frac{1}{p^{\prime }}=1-\frac{1%
}{p}\right) \text{, where }\psi ^{\varepsilon }\left( x\right) = \\ 
\psi \left( x,\frac{x}{\varepsilon }\right) \text{ }\left( x\in \Omega
\right) \text{;}%
\end{array}%
\end{equation*}

(ii) strongly $\Sigma $-converge in $L^{p}\left( \Omega \right) $ to some $%
u_{0}\in L^{p}\left( \Omega ;L_{per}^{p}\left( Y\right) \right) $ if the
following property is verified: 
\begin{equation*}
\left\{ 
\begin{array}{c}
\text{Given }\eta >0\text{ and }v\in L^{p}\left( \Omega ;\mathcal{C}%
_{per}\left( Y\right) \right) \text{ with} \\ 
\left\Vert u_{0}-v\right\Vert _{L^{p}\left( \Omega \times Y\right) }\leq 
\frac{\eta }{2}\text{, there is some }\alpha >0\text{ such} \\ 
\text{that }\left\Vert u_{\varepsilon }-v^{\varepsilon }\right\Vert
_{L^{p}\left( \Omega \right) }\leq \eta \text{ provided }E\ni \varepsilon
\leq \alpha \text{.}%
\end{array}%
\right.
\end{equation*}
\end{definition}

We will briefly express weak and strong $\Sigma $-convergence by writing $%
u_{\varepsilon }\rightarrow u_{0}$ in $L^{p}\left( \Omega \right) $-weak $%
\Sigma $ and $u_{\varepsilon }\rightarrow u_{0}$ in $L^{p}\left( \Omega
\right) $-strong $\Sigma $, respectively. Instead of repeating here the main
results underlying $\Sigma $-convergence theory for periodic structures, we
find it more convenient to draw the reader's attention to a few references
regarding two-scale convergence, e.g., \cite{bib1}, \cite{bib2}, \cite{bib7}
and \cite{bib8}. However, we recall below two fundamental results.

\begin{theorem}
\label{th2.1} Assume that $1<p<\infty $ and further $E$ is a fundamental
sequence. Let a sequence $\left( u_{\varepsilon }\right) _{\varepsilon \in
E} $ be bounded in $L^{p}\left( \Omega \right) $. Then, a subsequence $%
E^{\prime }$ can be extracted from $E$ such that $\left( u_{\varepsilon
}\right) _{\varepsilon \in E^{\prime }}$ weakly $\Sigma $-converges in $%
L^{p}\left( \Omega \right) $.
\end{theorem}

\begin{theorem}
\label{th2.2} Let $E$ be a fundamental sequence. Suppose a sequence $\left(
u_{\varepsilon }\right) _{\varepsilon \in E}$ is bounded in $H^{1}\left(
\Omega \right) =W^{1,2}\left( \Omega \right) $. Then, a subsequence $%
E^{\prime }$ can be extracted from $E$ such that, as $E^{\prime }\ni
\varepsilon \rightarrow 0$, 
\begin{equation*}
u_{\varepsilon }\rightarrow u_{0}\text{ in }H^{1}\left( \Omega \right) \text{%
-weak,\qquad \qquad \qquad \qquad \qquad \qquad }
\end{equation*}%
\begin{equation*}
u_{\varepsilon }\rightarrow u_{0}\text{ in }L^{2}\left( \Omega \right) \text{%
-weak }\Sigma \text{,\qquad \qquad \qquad \qquad \qquad }\quad
\end{equation*}%
\begin{equation*}
\frac{\partial u_{\varepsilon }}{\partial x_{j}}\rightarrow \frac{\partial
u_{0}}{\partial x_{j}}+\frac{\partial u_{1}}{\partial y_{j}}\text{ in }%
L^{2}\left( \Omega \right) \text{-weak }\Sigma \text{ }\left( 1\leq j\leq
N\right) \text{,}
\end{equation*}%
where $u_{0}\in H^{1}\left( \Omega \right) $, $u_{1}\in L^{2}\left( \Omega
;H_{\#}^{1}\left( Y\right) \right) $.
\end{theorem}

Now, let us also introduce the notion of \textit{two-scale convergence on
periodic surfaces. We denote by }$L_{per}^{p}\left( \partial T\right) $ the
space of functions $u$ in $L_{loc}^{p}\left( \partial \Theta \right) $
verifying $u\left( y+k\right) =u\left( y\right) $ for all $k\in \mathbb{Z}%
^{N}$ and for almost all $y\in \partial \Theta $ ($\partial \Theta $ is the
boundary of $\Theta $). Let $\partial T^{\varepsilon }$ be the boundary of $%
T^{\varepsilon }$ ($T^{\varepsilon }$ is given by (\ref{eq1.1b})).

\begin{definition}
\label{def2.2} A sequence $\left( u_{\varepsilon }\right) _{\varepsilon \in
E}$ with $u_{\varepsilon }\in L^{p}\left( \partial T^{\varepsilon }\right) $
for all $\varepsilon \in E$ is said to two-scale converge to some $u_{0}\in
L^{p}\left( \Omega ;L_{per}^{p}\left( \partial T\right) \right) $ if as $%
E\ni \varepsilon \rightarrow 0$,%
\begin{equation}
\varepsilon \int_{\partial T^{\varepsilon }}u_{\varepsilon }\left( x\right)
\psi ^{\varepsilon }\left( x\right) d\sigma _{\varepsilon }\left( x\right)
\rightarrow \int \int_{\Omega \times \partial T}u_{0}\left( x,y\right) \psi
\left( x,y\right) dxd\sigma \left( y\right)  \label{eq2.9a}
\end{equation}%
\begin{equation*}
\begin{array}{c}
\text{for all }\psi \in \mathcal{C}\left( \overline{\Omega };\mathcal{C}%
_{per}\left( Y\right) \right) \text{, where }\psi ^{\varepsilon }\left(
x\right) = \\ 
\psi \left( x,\frac{x}{\varepsilon }\right) \text{ }\left( x\in \Omega
\right)%
\end{array}%
\end{equation*}%
and where $d\sigma _{\varepsilon }$ and $d\sigma $ denote the surface
measures on $\partial T^{\varepsilon }$ and $\partial T$, respectively.
\end{definition}

The following result of convergence on periodic surfaces holds true.

\begin{theorem}
\label{th2.3} Let $1<p<+\infty $, and let $\left( u_{\varepsilon }\right)
_{\varepsilon \in E}$ be a sequence with $u_{\varepsilon }\in L^{p}\left(
\partial T^{\varepsilon }\right) $ for all $\varepsilon \in E$. Suppose that 
\begin{equation}
\varepsilon \int_{\partial T^{\varepsilon }}\left\vert u_{\varepsilon
}\left( x\right) \right\vert ^{p}d\sigma _{\varepsilon }\left( x\right) \leq
C  \label{eq2.9b}
\end{equation}
for all $\varepsilon \in E$, where $C$ is a constant independent of $%
\varepsilon $. Then, there exists a subsequence $E^{\prime }$ extracted from 
$E$ and a function $u_{0}\in L^{p}\left( \Omega ;L_{per}^{p}\left( \partial
T\right) \right) $ such that $\left( u_{\varepsilon }\right) _{\varepsilon
\in E^{\prime }}$ two-scale converges to $u_{0}$.
\end{theorem}

The proof of the preceding theorem can be found in, e.g. \cite{bib12} (see
also \cite{bib2}).

\begin{remark}
\label{rem2.2} It is of interest to know that if $u_{\varepsilon
}\rightarrow u_{0}$ in $L^{p}\left( \Omega \right) $-weak $\Sigma $, then (%
\ref{eq2.9}) holds for $\psi \in \mathcal{C}\left( \overline{\Omega }%
;L_{per}^{\infty }\left( Y\right) \right) $ (see \cite[Proposition 10]{bib9}
for the proof). Moreover if $\left( u_{\varepsilon }\right) _{\varepsilon
>0} $ with $u_{\varepsilon }\in L^{p}\left( \partial T^{\varepsilon }\right) 
$ two-scale converges to $u_{0}\in L^{p}\left( \Omega ;L_{per}^{p}\left(
\partial T\right) \right) $ (in the sense of Definition \ref{def2.2}), then (%
\ref{eq2.9a}) holds for $\psi \in \mathcal{C}\left( \overline{\Omega }%
;L_{per}^{\infty }\left( Y\right) \right) $. The reader can refer to \cite%
{bib12} for more details.
\end{remark}

The following useful proposition has its proof in \cite{bib2}.

\begin{proposition}
\label{pr2.2} Let $\left( u_{\varepsilon }\right) _{\varepsilon \in E}$ be a
sequence in $H^{1}\left( \Omega \right) $ such that 
\begin{equation*}
\left\Vert u_{\varepsilon }\right\Vert _{L^{2}\left( \Omega \right)
}+\varepsilon \left\Vert \nabla u_{\varepsilon }\right\Vert _{L^{2}\left(
\Omega \right) }\leq C\text{,}
\end{equation*}%
where $C>0$ is a constant independent of $\varepsilon $. Then the trace of $%
u_{\varepsilon }$ on $\partial T^{\varepsilon }$ satisfies%
\begin{equation*}
\varepsilon \int_{\partial T^{\varepsilon }}\left\vert u_{\varepsilon
}\left( x\right) \right\vert ^{2}d\sigma _{\varepsilon }\left( x\right) \leq
C
\end{equation*}%
for all $\varepsilon \in E$, and up to a subsequence, it two-scale converges
in the sens of Definition \ref{def2.2} to some $u_{0}\in L^{2}\left( \Omega
;L_{per}^{2}\left( \partial T\right) \right) $, which is the trace on $%
\partial T$ of a function in $L^{2}\left( \Omega ;H_{\#}^{1}\left( Y\right)
\right) $. More precisely, there exists a subsequence $E^{\prime }$ of $E$
and a function $u_{0}\in L^{2}\left( \Omega ;H_{per}^{1}\left( Y\right)
\right) $ such that as $E^{\prime }\ni \varepsilon \rightarrow 0$,%
\begin{equation*}
\varepsilon \int_{\partial T^{\varepsilon }}u_{\varepsilon }\left( x\right)
\psi ^{\varepsilon }\left( x\right) d\sigma _{\varepsilon }\left( x\right)
\rightarrow \int \int_{\Omega \times \partial T}u_{0}\left( x,y\right) \psi
\left( x,y\right) dxd\sigma \left( y\right) \text{ \ for all }\psi \in 
\mathcal{C}\left( \overline{\Omega };\mathcal{C}_{per}\left( Y\right)
\right) \text{,}
\end{equation*}%
\begin{equation*}
\int_{\Omega }u_{\varepsilon }\left( x\right) \psi ^{\varepsilon }\left(
x\right) dx\rightarrow \int \int_{\Omega \times Y}u_{0}\left( x,y\right)
\psi \left( x,y\right) dxdy\text{ \ for all }\psi \in L^{2}\left( \Omega ;%
\mathcal{C}_{per}\left( Y\right) \right)
\end{equation*}%
and 
\begin{equation*}
\varepsilon \int_{\Omega }\frac{\partial u_{\varepsilon }}{\partial x_{j}}%
\left( x\right) \psi ^{\varepsilon }\left( x\right) dx\rightarrow \int
\int_{\Omega \times Y}\frac{\partial u_{0}}{\partial y_{j}}\left( x,y\right)
\psi \left( x,y\right) dxdy\text{ \ for all }\psi \in L^{2}\left( \Omega ;%
\mathcal{C}_{per}\left( Y\right) \right) \text{,}
\end{equation*}%
for all $1\leq j\leq N$.
\end{proposition}

Having made the above preliminaries, let us turn now to the statement of the
hypotheses for the homogenization problem of (\ref{eq1.3})-(\ref{eq1.7}). In
view of (\ref{eq1.8a}) and since the functions $a_{ijkh}$, $\mathcal{\theta }
$ and $f_{j}$ belong to $L^{\infty }\left( \mathbb{R}^{N}\right) $ we have 
\begin{equation}
a_{ijkh}\text{, }\mathcal{\theta }\text{ and }f_{j}\in L_{per}^{\infty
}\left( Y\right) \text{\qquad }\left( 1\leq i,j,k,h\leq N\right) \text{.}
\label{eq2.10c}
\end{equation}%
Further, since the sets $k+\overset{\circ }{T}$ $\left( k\in \mathbb{Z}^{N}%
\text{, }\overset{\circ }{T}\text{ being the interior of }T\right) $ are
pairwise disjoint, the characteristic function, $\mathcal{\chi }_{\Theta }$,
of the set $\Theta $ ($\Theta $ is defined in (\ref{eq2.6a})) verifies 
\begin{equation*}
\mathcal{\chi }_{\Theta }=\sum_{k\in \mathbb{Z}^{N}}\mathcal{\chi }_{k+T}%
\text{ \quad a.e. in }\mathbb{R}^{N}\text{,}
\end{equation*}%
where $\mathcal{\chi }_{k+T}$ is the characteristic function of $k+T$ in $%
\mathbb{R}_{y}^{N}$. We have the following proposition.

\begin{proposition}
\label{pr2.3} The characteristic function of the set $\Theta $ ($\Theta $ is
given by (\ref{eq2.6a})), $\mathcal{\chi }_{\Theta }$ belongs to $%
L_{per}^{\infty }\left( Y\right) $ and moreover its mean value is 
\begin{equation*}
\int_{Y}\mathcal{\chi }_{\Theta }\left( y\right) dy=\left\vert T\right\vert 
\text{.}
\end{equation*}
\end{proposition}

Now, let%
\begin{equation*}
G=\mathbb{R}_{y}^{N}\backslash \Theta
\end{equation*}%
and $\mathcal{\chi }_{G}$ its characteristic function. We have:%
\begin{equation}
\mathcal{\chi }_{G}\in L_{per}^{\infty }\left( Y\right) \text{.}
\label{eq2.10a}
\end{equation}%
Indeed, $\mathcal{\chi }_{G}=1-\mathcal{\chi }_{\Theta }$ and $\mathcal{\chi 
}_{\Theta }\in L_{per}^{\infty }\left( Y\right) $ (in view of Proposition %
\ref{pr2.3}). Further, as $\int_{Y}\mathcal{\chi }_{\Theta }\left( y\right)
dy=\left\vert T\right\vert $ we have 
\begin{equation}
\int_{Y}\mathcal{\chi }_{G}\left( y\right) dy=\left\vert Y\right\vert
-\left\vert T\right\vert =\left\vert Y^{\ast }\right\vert \text{.}
\label{eq2.10b}
\end{equation}%
Moreover, let us notice that 
\begin{equation*}
\mathcal{\chi }_{\Theta }=\sum_{k\in \mathbb{Z}^{N}}\mathcal{\chi }_{k+%
\overset{\circ }{T}}
\end{equation*}%
a.e. in $\mathbb{R}^{N}$ and therefore $\mathcal{\chi }_{Y}\mathcal{\chi }%
_{\Theta }=\mathcal{\chi }_{\overset{\circ }{T}}=\mathcal{\chi }_{T}$ a.e.
in $\mathbb{R}^{N}$. Thus, we have 
\begin{equation}
\mathcal{\chi }_{Y}\mathcal{\chi }_{G}=\mathcal{\chi }_{Y^{\ast }}
\label{eq2.10d}
\end{equation}%
a.e. in $\mathbb{R}^{N}$.

The following useful lemma is proved in \cite{bib12}.

\begin{lemma}
\label{lem2.4} Let $E$ be a fundamental sequence. Let $\left( u_{\varepsilon
}\right) _{\varepsilon \in E}\subset L^{2}\left( \Omega \right) $ and $%
\left( v_{\varepsilon }\right) _{\varepsilon \in E}\subset L^{\infty }\left(
\Omega \right) $ be two sequences such that:

(i) $u_{\varepsilon }\rightarrow u_{0}$ in $L^{2}\left( \Omega \right) $%
-weak $\Sigma $ as $E\ni \varepsilon \rightarrow 0$,

(ii) $v_{\varepsilon }\rightarrow v_{0}$ in $L^{2}\left( \Omega \right) $%
-strong $\Sigma $ as $E\ni \varepsilon \rightarrow 0$,

(iii) $\left( v_{\varepsilon }\right) _{\varepsilon \in E}$ is bounded in $%
L^{\infty }\left( \Omega \right) $.

Then $u_{\varepsilon }v_{\varepsilon }\rightarrow u_{0}v_{0}$ in $%
L^{2}\left( \Omega \right) $-weak $\Sigma $ as $E\ni \varepsilon \rightarrow
0$.
\end{lemma}

\section{A convergence result for the homogenization process}

In the present section, our goal is to investigate the limiting behaviour,
as $\varepsilon \rightarrow 0$, of $\mathbf{u}_{\varepsilon }$ solution to (%
\ref{eq1.3})-(\ref{eq1.7}). To this end, let us state some preliminaries.

We have the following proposition on the estimates of solutions to (\ref%
{eq1.3})-(\ref{eq1.7}).

\begin{proposition}
\label{pr3.1} Suppose that (\ref{eq1.2a})-(\ref{eq1.2}) are verified. For $%
0<\varepsilon <1$, let $\mathbf{u}_{\varepsilon }$ be the unique solution to
(\ref{eq1.3})-(\ref{eq1.7}) and $\mathcal{P}^{\varepsilon }$ the extension
operator in Proposition \ref{pr2.1}. There exists a constant $C>0$
independent of $\varepsilon $ such that 
\begin{equation}
\left\Vert \mathcal{P}^{\varepsilon }\mathbf{u}_{\varepsilon }\right\Vert
_{H^{1}\left( \Omega \right) ^{N}}\leq C\text{,}  \label{eq2.12a}
\end{equation}%
\begin{equation}
\varepsilon \sum_{k=1}^{N}\int_{\partial T^{\varepsilon }}\left\vert
u_{\varepsilon }^{k}\left( x\right) \right\vert ^{2}d\sigma _{\varepsilon
}\left( x\right) \leq C\text{.}  \label{eq2.12b}
\end{equation}
\end{proposition}

\begin{proof}
We turn back to the variational probem (\ref{eq1.8}), we take in particular $%
\mathbf{v}=\mathbf{u}_{\varepsilon }$ in the equation. In virtue of (\ref%
{eq1.2a})-(\ref{eq1.2}), this leads to 
\begin{equation}
\alpha \left\Vert \mathbf{e}\left( \mathbf{u}_{\varepsilon }\right)
\right\Vert _{L^{2}\left( \Omega ^{\varepsilon }\right) }^{2}+\varepsilon
\alpha _{0}\sum_{k=1}^{N}\int_{\partial T^{\varepsilon }}\left\vert
u_{\varepsilon }^{k}\left( x\right) \right\vert ^{2}d\sigma _{\varepsilon
}\left( x\right) \leq \left\Vert \mathbf{u}_{\varepsilon }\right\Vert
_{L^{2}\left( \Omega ^{\varepsilon }\right) }\left\Vert \mathbf{f}%
^{\varepsilon }\right\Vert _{L^{2}\left( \Omega \right) }\text{.}
\label{eq2.13a}
\end{equation}%
But, by Proposition \ref{pr2.1}, (\ref{eq2.13a}) leads to 
\begin{equation*}
\frac{\alpha }{c^{\prime 2}c^{2}}\left\Vert \mathcal{P}^{\varepsilon }%
\mathbf{u}_{\varepsilon }\right\Vert _{L^{2}\left( \Omega \right) }^{2}\leq
\lambda \left( \Omega \right) ^{\frac{1}{2}}\left\Vert \mathbf{f}\right\Vert
_{\infty }\left\Vert \mathcal{P}^{\varepsilon }\mathbf{u}_{\varepsilon
}\right\Vert _{L^{2}\left( \Omega \right) }\text{,}
\end{equation*}%
where $c$ is the constant in (\ref{eq2.3b}) and $c^{\prime }$ the one in the
Korn's inequality. It follows by (\ref{eq2.13a}), Proposition \ref{pr2.1}
and the preceding inequality that%
\begin{equation}
\left\Vert \mathcal{P}^{\varepsilon }\mathbf{u}_{\varepsilon }\right\Vert
_{L^{2}\left( \Omega \right) }\leq \frac{c^{\prime 2}c^{2}}{\alpha }\lambda
\left( \Omega \right) ^{\frac{1}{2}}\left\Vert \mathbf{f}\right\Vert
_{\infty }\text{, }\left\Vert \mathbf{e}\left( \mathcal{P}^{\varepsilon }%
\mathbf{u}_{\varepsilon }\right) \right\Vert _{L^{2}\left( \Omega \right)
}\leq \frac{c^{2}c^{\prime }}{\alpha }\lambda \left( \Omega \right)
\left\Vert \mathbf{f}\right\Vert _{\infty }\text{ }  \label{eq2.13d}
\end{equation}%
and 
\begin{equation}
\varepsilon \sum_{k=1}^{N}\int_{\partial T^{\varepsilon }}\left\vert
u_{\varepsilon }^{k}\left( x\right) \right\vert ^{2}d\sigma _{\varepsilon
}\left( x\right) \leq \frac{c^{\prime 2}c^{2}}{\alpha _{0}\alpha }\lambda
\left( \Omega \right) \left\Vert \mathbf{f}\right\Vert _{\infty }^{2}
\label{eq2.13e}
\end{equation}%
for all $0<\varepsilon <1$. In view of (\ref{eq2.13d})-(\ref{eq2.13e}),
there exists a positive constant $C$ independent of $\varepsilon $ such that
(\ref{eq2.12a})-(\ref{eq2.12b}) are satisfied.
\end{proof}

Before we can establish the so-called global homogenization theorem for (\ref%
{eq1.3})-(\ref{eq1.7}), we require a few basic notation and results. To
begin, let 
\begin{equation*}
\mathcal{V}=\left\{ \mathbf{v}\in \mathcal{D}\left( \overline{\Omega };%
\mathbb{R}\right) ^{N}:\mathbf{v}=0\text{ on }\Gamma _{1}\right\}
\end{equation*}%
($\mathcal{D}\left( \overline{\Omega };\mathbb{R}\right) $ being the
restrictions to $\Omega $ of functions in $\mathcal{D}\left( \mathbb{R}^{N};%
\mathbb{R}\right) $),%
\begin{equation*}
\mathbf{V}_{0}=\left\{ \mathbf{v}\in H^{1}\left( \Omega ;\mathbb{R}\right)
^{N}:\mathbf{v}=0\text{ on }\Gamma _{1}\right\} \text{,}
\end{equation*}%
\begin{equation*}
\mathcal{V}_{Y}=\left\{ \mathbf{\psi }\in \mathcal{C}_{per}^{\infty }\left(
Y;\mathbb{R}\right) ^{N}:\int_{Y}\mathbf{\psi }\left( y\right) dy=0\right\}
\end{equation*}%
and 
\begin{equation*}
\mathbf{V}_{Y}=H_{\#}^{1}\left( Y;\mathbb{R}\right) ^{N}\text{ }
\end{equation*}%
where: $\mathcal{C}_{per}^{\infty }\left( Y;\mathbb{R}\right) =\mathcal{C}%
^{\infty }\left( \mathbb{R}^{N};\mathbb{R}\right) \cap \mathcal{C}%
_{per}\left( Y\right) $. We provide $\mathbf{V}_{Y}$ with the $%
H_{\#}^{1}\left( Y\right) ^{N}$-norm, which makes it a Hilbert space. There
is no difficulty in verifying that $\mathcal{V}_{Y}$ is dense in $\mathbf{V}%
_{Y}$. With this in mind, we denote by $\mathbf{V}$ de closure of $\mathcal{V%
}$ in $H^{1}\left( \Omega ;\mathbb{R}\right) ^{N}$ and we set 
\begin{equation*}
\mathbb{F}_{0}^{1}=\mathbf{V}\times L^{2}\left( \Omega ;\mathbf{V}%
_{Y}\right) \text{.}
\end{equation*}%
This is a Hilbert space with norm 
\begin{equation*}
\left\Vert \mathbf{v}\right\Vert _{\mathbb{F}_{0}^{1}}=\left( \left\Vert 
\mathbf{e}\left( \mathbf{v}_{0}\right) \right\Vert _{L^{2}\left( \Omega
\right) ^{N^{2}}}^{2}+\left\Vert \mathbf{e}^{y}\left( \mathbf{v}_{1}\right)
\right\Vert _{L^{2}\left( \Omega \times Y\right) ^{N^{2}}}^{2}\right) ^{%
\frac{1}{2}}\text{, }\mathbf{v=}\left( \mathbf{v}_{0},\mathbf{v}_{1}\right)
\in \mathbb{F}_{0}^{1}
\end{equation*}%
where $\mathbf{e}^{y}\left( \mathbf{v}_{1}\right) =\left( \mathbf{e}%
_{ij}^{y}\left( \mathbf{v}_{1}\right) \right) _{1\leq i,j\leq N}$ with%
\begin{equation*}
\mathbf{e}_{ij}^{y}\left( \mathbf{v}_{1}\right) =\frac{1}{2}\left( \frac{%
\partial v_{1}^{i}}{\partial y_{j}}+\frac{\partial v_{1}^{j}}{\partial y_{i}}%
\right) \text{, \qquad }\left( 1\leq i,j\leq N\right) \text{.}
\end{equation*}%
On the other hand, we put 
\begin{equation*}
\mathbf{\tciFourier }_{0}^{\infty }=\mathcal{V\times }\left[ \mathcal{D}%
\left( \Omega ;\mathbb{R}\right) \otimes \mathcal{V}_{Y}\right] \text{,}
\end{equation*}%
where $\mathcal{D}\left( \Omega ;\mathbb{R}\right) \otimes \mathcal{V}_{Y}$
stands for the space of vector functions $\mathbf{\phi }$ on $\Omega \times 
\mathbb{R}_{y}^{N}$ of the form 
\begin{equation*}
\mathbf{\phi }\left( x,y\right) =\sum_{finite}\varphi _{i}\left( x\right) 
\mathbf{w}_{i}\left( y\right) \text{ }\left( x\in \Omega ,\text{ }y\in 
\mathbb{R}^{N}\right)
\end{equation*}%
with $\varphi _{i}\in \mathcal{D}\left( \Omega ;\mathbb{R}\right) $, $%
\mathbf{w}_{i}\in \mathcal{V}_{Y}$. It is clear that $\mathbf{\tciFourier }%
_{0}^{\infty }$ is dense in $\mathbb{F}_{0}^{1}$.

It is of interest to notice that for $\mathbf{v}=\left( \mathbf{v}_{0},%
\mathbf{v}_{1}\right) \in \mathbb{F}_{0}^{1}$ with $\mathbf{v}_{0}=\left(
v_{0}^{k}\right) _{1\leq k\leq N}$ and $\mathbf{v}_{1}=\left(
v_{1}^{k}\right) _{1\leq k\leq N}$, if we set

\begin{equation*}
\mathbf{E}_{ij}\left( \mathbf{v}\right) =\mathbf{e}_{ij}\left( \mathbf{v}%
_{0}\right) +\mathbf{e}_{ij}^{y}\left( \mathbf{v}_{1}\right) \text{ }\qquad
\left( 1\leq i,j\leq N\right) \text{,}
\end{equation*}%
then we have 
\begin{equation}
\left\Vert \mathbf{v}\right\Vert _{\mathbb{F}_{0}^{1}}=\left(
\sum_{i,j=1}^{N}\left\Vert \mathbf{E}_{ij}\left( \mathbf{v}\right)
\right\Vert _{L^{2}\left( \Omega \times Y\right) }^{2}\right) ^{\frac{1}{2}%
}\qquad \left( \mathbf{v}=\left( v_{0},v_{1}\right) \in \mathbb{F}%
_{0}^{1}\right) \text{.}  \label{eq2.11}
\end{equation}

Now, for $\mathbf{u}=\left( \mathbf{u}_{0},\mathbf{u}_{1}\right) $ and $%
\mathbf{v}=\left( \mathbf{v}_{0},\mathbf{v}_{1}\right) \in \mathbb{F}%
_{0}^{1} $ we set 
\begin{equation*}
\widehat{a}_{\Omega }\left( \mathbf{u},\mathbf{v}\right)
=\sum_{i,j,k,h=1}^{N}\int \int_{\Omega \times Y^{\ast }}a_{ijkh}\mathbf{E}%
_{ij}\left( \mathbf{u}\right) \mathbf{E}_{kh}\left( \mathbf{v}\right)
dxdy+\int \int_{\Omega \times \partial T}\mathcal{\theta }\mathbf{u}_{0}%
{\small \cdot }\mathbf{v}_{0}dxd\sigma \text{.}
\end{equation*}%
This defines a bilinear form $\widehat{a}_{\Omega }\left( ,\right) $ on $%
\mathbb{F}_{0}^{1}\times \mathbb{F}_{0}^{1}$ which in view of (\ref{eq1.2a}%
)-(\ref{eq1.2}), is symmetric, positive, continuous and noncoercive. Indeed,
for some $\mathbf{u}=\left( \mathbf{u}_{0},\mathbf{u}_{1}\right) \in \mathbb{%
F}_{0}^{1}$, $\widehat{a}_{\Omega }\left( \mathbf{u},\mathbf{u}\right) =0$
if and only if $\mathbf{u}_{0}=0$ and $\mathbf{e}_{ij}^{y}\left( \mathbf{v}%
_{1}\right) \left( x,y\right) =0$ a.e. in $\left( x,y\right) \in \Omega
\times Y^{\ast }$ ($1\leq i,j\leq N$). Thus, $\mathbf{u}=\left( \mathbf{u}%
_{0},\mathbf{u}_{1}\right) $ is not necessarily the zero function in $%
\mathbb{F}_{0}^{1}$. However, we put 
\begin{equation*}
\mathbf{N}\left( \mathbf{v}\right) =\left( \sum_{i,j=1}^{N}\int \int_{\Omega
\times Y^{\ast }}\left\vert \mathbf{E}_{ij}\left( \mathbf{v}\right)
\right\vert ^{2}dxdy+\int_{\Omega }\left\vert \mathbf{v}_{0}\right\vert
^{2}dx\right) ^{\frac{1}{2}}
\end{equation*}%
for all $\mathbf{v}=\left( \mathbf{v}_{0},\mathbf{v}_{1}\right) \in \mathbb{F%
}_{0}^{1}$. This defines a seminorm on $\mathbb{F}_{0}^{1}$. Equipped with
the seminorm $\mathbf{N}\left( .\right) $, $\mathbb{F}_{0}^{1}$ is a
pre-Hilbert space which is nonseparated and noncomplete. Further, let us
consider the linear form $\widehat{l}_{\Omega }$ on $\mathbb{F}_{0}^{1}$
defined by 
\begin{equation*}
\widehat{l}_{\Omega }\left( \mathbf{v}\right) =\int \int_{\Omega \times
Y^{\ast }}\mathbf{f}\left( y\right) {\small \cdot }\mathbf{v}_{0}\left(
x\right) dxdy+\int_{\Gamma _{2}}\mathbf{t}{\small \cdot }\mathbf{v}%
_{0}d\Gamma
\end{equation*}%
for all $\mathbf{v}=\left( \mathbf{v}_{0},\mathbf{v}_{1}\right) \in \mathbb{F%
}_{0}^{1}$. The form $\widehat{l}_{\Omega }$ is continuous on $\mathbb{F}%
_{0}^{1}$ for the norm $\left\Vert .\right\Vert _{\mathbb{F}_{0}^{1}}$ and
the seminorm $\mathbf{N}\left( .\right) $. Therefore, we have the following
lemma.

\begin{lemma}
\label{lem3.2} Suppose (\ref{eq1.2a})-(\ref{eq1.2}) and (\ref{eq1.8a}) hold.
There exists $\mathbf{u=}\left( \mathbf{u}_{0},\mathbf{u}_{1}\right) \in 
\mathbb{F}_{0}^{1}$ satisfying the variational problem%
\begin{equation}
\widehat{a}_{\Omega }\left( \mathbf{u},\mathbf{v}\right) =\widehat{l}%
_{\Omega }\left( \mathbf{v}\right) \text{ \quad for all }\mathbf{v}\in 
\mathbb{F}_{0}^{1}\text{.}  \label{eq2.12}
\end{equation}%
Moreover, $\mathbf{u}_{0}$ is strictly unique and $\mathbf{u}_{1}$ is unique
up to an additive vector function $\mathbf{g}\in L^{2}\left( \Omega
;V_{Y}\right) $ such that $\mathbf{e}_{ij}^{y}\left( \mathbf{g}\right)
\left( x,y\right) =0$ a.e. in $\left( x,y\right) \in \Omega \times Y^{\ast }$
$\left( 1\leq i,j\leq N\right) $.
\end{lemma}

\begin{proof}
The proof of this lemma is a simple adaptation of the one in \cite[Lemma 2.5]%
{bib10}. So, for shortness we omit it.
\end{proof}

Now, let us state our convergence theorem.

\begin{theorem}
\label{th3.1} Suppose that the hypotheses (\ref{eq1.2a})-(\ref{eq1.2}) and (%
\ref{eq1.8a}) are satisfied. For $\varepsilon \in E$, let $\mathbf{u}%
_{\varepsilon }\in \mathbf{V}_{\varepsilon }$ be the unique solution to (\ref%
{eq1.3})-(\ref{eq1.7}) ($E$ being a fundamental sequence), and let $\mathcal{%
P}^{\varepsilon }$ be the extension operator of Proposition \ref{pr2.1}.
Then, a subsequence $E^{\prime }$ can be extracted from $E$ such that as $%
E^{\prime }\ni \varepsilon \rightarrow 0$, 
\begin{equation}
\mathcal{P}^{\varepsilon }\mathbf{u}_{\varepsilon }\rightarrow \mathbf{u}_{0}%
\text{ in }H^{1}\left( \Omega \right) ^{N}\text{-weak,}  \label{eq2.15a}
\end{equation}%
\begin{equation}
\mathbf{e}_{ij}\left( \mathcal{P}^{\varepsilon }\mathbf{u}_{\varepsilon
}\right) \rightarrow \mathbf{e}_{ij}\left( \mathbf{u}_{0}\right) +\mathbf{e}%
_{ij}^{y}\left( \mathbf{u}_{1}\right) \text{ in }L^{2}\left( \Omega \right) 
\text{-weak }\Sigma \text{, \quad }\left( 1\leq i,j\leq N\right)
\label{eq2.15c}
\end{equation}%
where $\mathbf{u}_{1}\in L^{2}\left( \Omega ;\mathbf{V}_{Y}\right) $ and $%
\mathbf{u}=\left( \mathbf{u}_{0},\mathbf{u}_{1}\right) $ verifies the
variational equality (\ref{eq2.12}).
\end{theorem}

\begin{proof}
Let $E$ be a fundamental sequence. According to Proposition \ref{pr3.1}, the
sequence $\left( \mathcal{P}^{\varepsilon }\mathbf{u}_{\varepsilon }\right)
_{\varepsilon \in E}$ is bounded in $H^{1}\left( \Omega \right) ^{N}$ in
view of (\ref{eq2.12a}). Thus, Theorems \ref{th2.1} and \ref{th2.2} yield a
subsequence $E^{\prime }$ extracted from $E$, a vector function $\mathbf{u}%
=\left( \mathbf{u}_{0},\mathbf{u}_{1}\right) \in H^{1}\left( \Omega ;\mathbb{%
R}\right) ^{N}\times L^{2}\left( \Omega ;H_{\#}^{1}\left( Y;\mathbb{R}%
\right) ^{N}\right) $ such that (\ref{eq2.15a})-(\ref{eq2.15c}) hold. Let us
check that $\mathbf{u}=\left( \mathbf{u}_{0},\mathbf{u}_{1}\right) $
verifies (\ref{eq2.12}). For each real $\varepsilon >0$, let 
\begin{equation*}
\mathbf{\phi }_{\varepsilon }=\mathbf{\phi }_{0}+\varepsilon \mathbf{\phi }%
_{1}^{\varepsilon }\text{ with }\mathbf{\phi }_{0}=\left( \phi
_{0}^{k}\right) _{1\leq k\leq N}\in \mathcal{V}\text{,\quad\ }\mathbf{\phi }%
_{1}=\left( \phi _{1}^{k}\right) _{1\leq k\leq N}\in \mathcal{D}\left(
\Omega ;\mathbb{R}\right) \otimes \mathcal{V}_{Y}\text{,}
\end{equation*}%
i.e., $\mathbf{\phi }_{\varepsilon }\left( x\right) =\mathbf{\phi }%
_{0}\left( x\right) +\varepsilon \mathbf{\phi }_{1}\left( x,\frac{x}{%
\varepsilon }\right) $ for $x\in \Omega $. Clearly, we have $\mathbf{\phi }%
_{\varepsilon }\in \mathcal{V}$. This being so, taking in (\ref{eq1.8}) $%
\mathbf{v}=\mathbf{\phi }_{\varepsilon }{\small \mid }_{\Omega ^{\varepsilon
}}$ ($\mathbf{\phi }_{\varepsilon }{\small \mid }_{\Omega ^{\varepsilon }}$
being de restriction of $\mathbf{\phi }_{\varepsilon }$ to $\Omega
^{\varepsilon }$) leads to 
\begin{equation}
\mathbf{a}_{\varepsilon }\left( \mathbf{u}_{\varepsilon },\mathbf{\phi }%
_{\varepsilon }{\small \mid }_{\Omega ^{\varepsilon }}\right) =\int_{\Omega
^{\varepsilon }}\mathbf{f}^{\varepsilon }{\small \cdot }\mathbf{\phi }%
_{\varepsilon }dx+\int_{\Gamma _{2}}\mathbf{t}{\small \cdot }\mathbf{\phi }%
_{\varepsilon }d\Gamma \text{.}  \label{eq2.13}
\end{equation}%
Further, by the decomposition $\Omega ^{\varepsilon }=Q^{\varepsilon }\cup
\left( \Omega ^{\varepsilon }{\small \diagdown }Q^{\varepsilon }\right) $
and use of $Q^{\varepsilon }=\Omega \cap \varepsilon G$, the equality (\ref%
{eq2.13}) yields 
\begin{equation}
\begin{array}{c}
\sum_{i,j,k,h=1}^{N}\int_{\Omega }a_{ijkh}^{\varepsilon }\mathbf{e}%
_{ij}\left( \mathcal{P}^{\varepsilon }\mathbf{u}_{\varepsilon }\right) 
\mathbf{e}_{kh}\left( \mathbf{\phi }_{\varepsilon }\right) \mathcal{\chi }%
_{G}^{\varepsilon }dx+\sum_{i,j,k,h=1}^{N}\int_{\Omega ^{\varepsilon }%
{\small \diagdown }Q^{\varepsilon }}a_{ijkh}^{\varepsilon }\mathbf{e}%
_{ij}\left( \mathcal{P}^{\varepsilon }\mathbf{u}_{\varepsilon }\right) 
\mathbf{e}_{kh}\left( \mathbf{\phi }_{\varepsilon }\right) dx \\ 
+\varepsilon \int_{\partial T^{\varepsilon }}\mathcal{\theta }^{\varepsilon }%
\mathcal{P}^{\varepsilon }\mathbf{u}_{\varepsilon }{\small \cdot }\mathbf{%
\phi }_{\varepsilon }d\sigma _{\varepsilon }=\int_{\Omega }\mathbf{f}%
^{\varepsilon }{\small \cdot }\mathbf{\phi }_{\varepsilon }\mathcal{\chi }%
_{G}^{\varepsilon }dx+\int_{\Omega ^{\varepsilon }{\small \diagdown }%
Q^{\varepsilon }}\mathbf{f}^{\varepsilon }{\small \cdot }\mathbf{\phi }%
_{\varepsilon }dx+\int_{\Gamma _{2}}\mathbf{t}{\small \cdot }\mathbf{\phi }%
_{0}d\Gamma \text{,}%
\end{array}
\label{eq2.14}
\end{equation}%
for all $0<\varepsilon <1$, since $\mathbf{\phi }_{1}^{\varepsilon }=0$ on $%
\Gamma $. On the other hand, in \cite[proof of Lemma 2.1]{bib12} it is
verified that 
\begin{equation*}
\Omega ^{\varepsilon }{\small \diagdown }Q^{\varepsilon }\subset
J^{\varepsilon }\left( \partial \Omega \right) \qquad \left( 0<\varepsilon
<1\right) \text{,}
\end{equation*}%
where $J^{\varepsilon }\left( \partial \Omega \right) =\cup _{k\in
j^{\varepsilon }\left( \partial \Omega \right) }\varepsilon \left( k+%
\overline{Y}\right) $ and where $j^{\varepsilon }\left( \partial \Omega
\right) =\left\{ k\in \mathbb{Z}^{N}:\varepsilon \left( k+\overline{Y}%
\right) \cap \partial \Omega \neq \varnothing \right\} $. Moreover, in
virtue of the regularity of the Lebesgue measure $\lambda $, we have $%
\lambda \left( J^{\varepsilon }\left( \partial \Omega \right) \right)
\rightarrow \lambda \left( \partial \Omega \right) =0$ as $\varepsilon
\rightarrow 0$. Thus as $\varepsilon \rightarrow 0$, $\lambda \left( \Omega
^{\varepsilon }{\small \diagdown }Q^{\varepsilon }\right) \rightarrow 0$ and
therefore

$\sum_{i,j,k,h=1}^{N}\int_{\Omega ^{\varepsilon }{\small \diagdown }%
Q^{\varepsilon }}a_{ijkh}^{\varepsilon }\mathbf{e}_{ij}\left( \mathcal{P}%
^{\varepsilon }\mathbf{u}_{\varepsilon }\right) \mathbf{e}_{kh}\left( 
\mathbf{\phi }_{\varepsilon }\right) dx$ and $\int_{\Omega ^{\varepsilon }%
{\small \diagdown }Q^{\varepsilon }}\mathbf{f}^{\varepsilon }{\small \cdot }%
\mathbf{\phi }_{\varepsilon }dx$ tend to $0$ as $\varepsilon \rightarrow 0$.
Let us pass to the limit in (\ref{eq2.14}) when $E^{\prime }\ni \varepsilon
\rightarrow 0$. First, by (\ref{eq2.10c}) and (\ref{eq2.10a}) we see that $%
a_{ijkh}\mathcal{\chi }_{G}$ and $f_{j}\mathcal{\chi }_{G}$ belong to $%
L_{per}^{\infty }\left( Y\right) $. Further, for $1\leq k,h\leq N$ the
sequence $\left( \mathbf{e}_{kh}\left( \mathbf{\phi }_{\varepsilon }\right)
\right) _{0<\varepsilon <1}$ is bounded in $L^{\infty }\left( \Omega \right) 
$ and 
\begin{equation*}
\mathbf{e}_{kh}\left( \mathbf{\phi }_{\varepsilon }\right) \rightarrow 
\mathbf{E}_{kh}\left( \mathbf{\phi }\right) =\mathbf{e}_{kh}\left( \mathbf{%
\phi }_{0}\right) +\mathbf{e}_{kh}^{y}\left( \mathbf{\phi }_{1}\right) \text{
in }L^{2}\left( \Omega \right) \text{-strong }\Sigma
\end{equation*}%
as $\varepsilon \rightarrow 0$ (see, e.g., \cite[Lemma 2.2]{bib9} for
details) with $\mathbf{\phi }=\left( \mathbf{\phi }_{0},\mathbf{\phi }%
_{1}\right) $. Then, according to (\ref{eq2.15c}) it follows by Lemma \ref%
{lem2.4} that 
\begin{equation*}
\mathbf{e}_{ij}\left( \mathcal{P}^{\varepsilon }\mathbf{u}_{\varepsilon
}\right) \mathbf{e}_{kh}\left( \mathbf{\phi }_{\varepsilon }\right)
\rightarrow \mathbf{E}_{ij}\left( \mathbf{u}\right) \mathbf{E}_{kh}\left( 
\mathbf{\phi }\right) \text{ in }L^{2}\left( \Omega \right) \text{-weak }%
\Sigma
\end{equation*}%
as $E^{\prime }\ni \varepsilon \rightarrow 0$ with $\mathbf{u}=\left( 
\mathbf{u}_{0},\mathbf{u}_{1}\right) $. Moreover, by Remark \ref{rem2.2} we
see that 
\begin{equation}
\sum_{i,j,k,h=1}^{N}\int_{\Omega }a_{ijkh}^{\varepsilon }\mathbf{e}%
_{ij}\left( \mathcal{P}^{\varepsilon }\mathbf{u}_{\varepsilon }\right) 
\mathbf{e}_{kh}\left( \mathbf{\phi }_{\varepsilon }\right) \mathcal{\chi }%
_{G}^{\varepsilon }dx\rightarrow \sum_{i,j,k,h=1}^{N}\int \int_{\Omega
\times Y}a_{ijkh}\mathcal{\chi }_{G}\mathbf{E}_{ij}\left( \mathbf{u}\right) 
\mathbf{E}_{kh}\left( \mathbf{\phi }\right) dxdy  \label{eq2.15}
\end{equation}%
as $E^{\prime }\ni \varepsilon \rightarrow 0$. On the other hand, by (\ref%
{eq2.12b}) of Proposition \ref{pr3.1}, Proposition \ref{pr2.2} and (\ref%
{eq2.15a}) the subsequence $E^{\prime }$ can be extracted from $E$ such that 
\begin{equation*}
\varepsilon \int_{\partial T^{\varepsilon }}\mathcal{\theta }^{\varepsilon }%
\mathcal{P}^{\varepsilon }\mathbf{u}_{\varepsilon }{\small \cdot }\mathbf{%
\phi }_{0}d\sigma _{\varepsilon }\rightarrow \int \int_{\Omega \times
\partial T}\mathcal{\theta }\left( y\right) \mathbf{u}_{0}\left( x\right) 
{\small \cdot }\mathbf{\phi }_{0}\left( x\right) dxd\sigma \left( y\right)
\end{equation*}%
as $E^{\prime }\ni \varepsilon \rightarrow 0$, in virtue of Remark \ref%
{rem2.2} since $\mathcal{\theta }\phi _{0}^{k}\in \mathcal{C}\left( 
\overline{\Omega };L_{per}^{\infty }\left( Y\right) \right) $. Furthermore,
by the same argument we see that 
\begin{equation*}
\varepsilon ^{2}\int_{\partial T^{\varepsilon }}\mathcal{\theta }%
^{\varepsilon }\mathcal{P}^{\varepsilon }\mathbf{u}_{\varepsilon }{\small %
\cdot }\mathbf{\phi }_{1}^{\varepsilon }d\sigma _{\varepsilon }\rightarrow 0
\end{equation*}%
as $E^{\prime }\ni \varepsilon \rightarrow 0$, since $\mathcal{\theta }\phi
_{1}^{k}\in \mathcal{C}\left( \overline{\Omega };L_{per}^{\infty }\left(
Y\right) \right) $. Thus, as $E^{\prime }\ni \varepsilon \rightarrow 0$, 
\begin{equation}
\varepsilon \int_{\partial T^{\varepsilon }}\mathcal{\theta }^{\varepsilon }%
\mathcal{P}^{\varepsilon }\mathbf{u}_{\varepsilon }{\small \cdot }\mathbf{%
\phi }_{\varepsilon }d\sigma _{\varepsilon }\rightarrow \int \int_{\Omega
\times \partial T}\mathcal{\theta }\left( y\right) \mathbf{u}_{0}\left(
x\right) {\small \cdot }\mathbf{\phi }_{0}\left( x\right) dxd\sigma \left(
y\right) \text{.}  \label{eq2.16}
\end{equation}%
Once more, we use Remark \ref{rem2.2} to have 
\begin{equation}
\int_{\Omega }\mathbf{f}^{\varepsilon }{\small \cdot }\mathbf{\phi }%
_{\varepsilon }\mathcal{\chi }_{G}^{\varepsilon }dx\rightarrow \int
\int_{\Omega \times Y}\mathbf{f}{\small \cdot }\mathbf{\phi }_{0}\mathcal{%
\chi }_{G}dxdy\text{, }  \label{eq2.17}
\end{equation}%
as $E^{\prime }\ni \varepsilon \rightarrow 0$. Finally, we pass to the limit
in (\ref{eq2.14}) as $E^{\prime }\ni \varepsilon \rightarrow 0$ and we
obtain by (\ref{eq2.15})-(\ref{eq2.17}) and (\ref{eq2.10d}), 
\begin{equation}
\widehat{a}_{\Omega }\left( \mathbf{u},\mathbf{\phi }\right) =\widehat{l}%
_{\Omega }\left( \mathbf{\phi }\right) \text{,}  \label{eq2.18}
\end{equation}%
for all $\mathbf{\phi }=\left( \mathbf{\phi }_{0},\mathbf{\phi }_{1}\right)
\in \mathbf{\tciFourier }_{0}^{\infty }$. Thus, using the density of $%
\mathbf{\tciFourier }_{0}^{\infty }$ in $\mathbb{F}_{0}^{1}$ and the
continuity of the forms $\widehat{a}_{\Omega }\left( ,\right) $ and $%
\widehat{l}_{\Omega }$, we see that $\mathbf{u}=\left( \mathbf{u}_{0},%
\mathbf{u}_{1}\right) $ verifies (\ref{eq2.12}). The proof of the theorem is
complete.
\end{proof}

Now, we introduce the bilinear form $\widehat{a}$ on $\mathbf{V}_{Y}\times 
\mathbf{V}_{Y}$ defined by%
\begin{equation*}
\widehat{a}\left( \mathbf{v},\mathbf{w}\right)
=\sum_{i,j,k,h=1}^{N}\int_{Y^{\ast }}a_{ijkh}\mathbf{e}_{ij}^{y}\left( 
\mathbf{v}\right) \mathbf{e}_{kh}^{y}\left( \mathbf{w}\right) dy
\end{equation*}%
for $\mathbf{v}=\left( v^{k}\right) $ and $\mathbf{w}=\left( w^{k}\right) $
in $\mathbf{V}_{Y}$, which is positive in view of (\ref{eq1.2}). Next, for
each couple of indices $1\leq i,j\leq N$, we consider the variational problem%
\begin{equation}
\left\{ 
\begin{array}{c}
\mathbf{\chi }^{ij}\in \mathbf{V}_{Y}:\qquad \qquad \qquad \qquad \\ 
\widehat{a}\left( \mathbf{\chi }^{ij},\mathbf{w}\right)
=\sum_{k,h=1}^{N}\int_{Y^{\ast }}a_{ijkh}\mathbf{e}_{kh}^{y}\left( \mathbf{w}%
\right) dy \\ 
\text{for all }\mathbf{w}=\left( w^{k}\right) \ \text{in }\mathbf{V}_{Y}%
\text{,\qquad }%
\end{array}%
\right.  \label{eq2.19}
\end{equation}%
which admits a solution $\mathbf{\chi }^{ij}$, unique up to an additive
vector function $\mathbf{g}=\left( g^{k}\right) \in \mathbf{V}_{Y}$ such
that $\mathbf{e}_{ij}^{y}\left( \mathbf{g}\right) =0$ a.e. in $Y^{\ast }$.

\begin{lemma}
\label{lem3.3} Under the hypotheses and notations of Theorem \ref{th3.1},
there is some vector function $\mathbf{g}=\left( g^{k}\right) \in
L^{2}\left( \Omega ;\mathbf{V}_{Y}\right) $ such that $\mathbf{e}%
_{ij}^{y}\left( \mathbf{g}\right) =0$ a.e. in $\Omega \times Y^{\ast }$ and 
\begin{equation}
\mathbf{u}_{1}\left( x,y\right) =-\sum_{i,j=1}^{N}\mathbf{e}_{ij}\left( 
\mathbf{u}_{0}\right) \mathbf{\chi }^{ij}\left( y\right) +\mathbf{g}\left(
x,y\right)  \label{eq2.20a}
\end{equation}%
almost everywhere in $\left( x,y\right) \in \Omega \times \mathbb{R}^{N}$.
\end{lemma}

\begin{proof}
In (\ref{eq2.12}), choose the test functions $\mathbf{v}=\left( \mathbf{v}%
_{0},\mathbf{v}_{1}\right) $ such that $\mathbf{v}_{0}=0$, $\mathbf{v}%
_{1}\left( x,y\right) =\varphi \left( x\right) \mathbf{w}\left( y\right) $
for $\left( x,y\right) \in \Omega \times \mathbb{R}^{N}$, where $\varphi \in 
\mathcal{D}\left( \Omega ;\mathbb{R}\right) $ and $\mathbf{w}\in \mathbf{V}%
_{Y}$. Then, almost everywhere in $x\in \Omega $, we have 
\begin{equation}
\left\{ 
\begin{array}{c}
\widehat{a}\left( \mathbf{u}_{1}\left( x,.\right) ,\mathbf{w}\right)
=-\sum_{i,j,k,h=1}^{N}\mathbf{e}_{ij}\left( \mathbf{u}_{0}\right)
\int_{Y^{\ast }}a_{ijkh}\mathbf{e}_{kh}^{y}\left( \mathbf{w}\right) dy \\ 
\text{for all }\mathbf{w}=\left( w^{k}\right) \in \mathbf{V}_{Y}\text{%
.\qquad \qquad \qquad \qquad \qquad \quad }\qquad \qquad%
\end{array}%
\right.  \label{eq2.20}
\end{equation}%
But it is clear that up to an additive function $\mathbf{g}\left( x\right)
=\left( g^{k}\left( x\right) \right) \in \mathbf{V}_{Y}$ such that $\mathbf{e%
}_{ij}^{y}\left( \mathbf{g}\right) =0$ a.e. in $Y^{\ast }$, $\mathbf{u}%
_{1}\left( x,.\right) $ (for fixed $x\in \Omega $) is the sole function in $%
\mathbf{V}_{Y}$ solving the variational equation (\ref{eq2.20}). On the
other hand, it is an easy matter to check that the function of $y$ on the
right of (\ref{eq2.20a}) solves the same variational problem. Hence the
lemma follows immediately.
\end{proof}

The next point deals with the so-called macroscopic homogenized equations
for (\ref{eq1.3})-(\ref{eq1.7}). Our goal here is to derive a well-posed
boundary value problem for $\mathbf{u}_{0}$. To begin, for $1\leq
i,j,k,h\leq N$, let 
\begin{equation}
q_{ijkh}=\int_{Y^{\ast }}a_{ijkh}\left( y\right)
dy-\sum_{p,q=1}^{N}\int_{Y^{\ast }}a_{pqkh}\left( y\right) \mathbf{e}%
_{pq}^{y}\left( \mathbf{\chi }^{ij}\right) \left( y\right) dy\text{.}
\label{eq2.21a}
\end{equation}%
where $\mathbf{\chi }^{ij}$ is given by (\ref{eq2.19}). There is no
difficulty to verify that 
\begin{equation*}
q_{ijkh}=\widehat{a}\left( \mathbf{\chi }^{ij}-\mathbf{\pi }^{ij},\mathbf{%
\chi }^{kh}-\mathbf{\pi }^{kh}\right)
\end{equation*}%
where $\mathbf{\pi }^{ij}=\left( \pi _{l}^{ij}\right) $ is the vector
function with $\pi _{l}^{ij}\left( y\right) =y_{i}\delta _{lj}$ for $1\leq
l\leq N$, $\delta _{lj}$ being the Kr\"{o}necker symbol. So the coefficients 
$q_{ijkh}$ have the following properties:%
\begin{equation}
q_{ijkh}=q_{khij}\text{,}  \label{eq2.21b}
\end{equation}%
\begin{equation}
\sum_{i,j,k,h=1}^{N}q_{ijkh}\zeta _{ij}\zeta _{kh}\geq 0  \label{eq2.21c}
\end{equation}%
for all $\zeta =\left( \zeta _{ij}\right) \in \mathbb{R}^{N\times N}$. Next,
we consider the following Hooke's type relation%
\begin{equation*}
\mathbf{\sigma }_{ij}^{0}=\sum_{k,h=1}^{N}q_{ijkh}\mathbf{e}_{kh}\left( 
\mathbf{u}_{0}\right) \text{,\qquad }\left( 1\leq i,j\leq N\right)
\end{equation*}%
and the following boundary value problem 
\begin{equation}
\qquad \qquad \mathcal{-}\func{div}\mathbf{\sigma }^{0}+\widetilde{\mathcal{%
\theta }}\mathbf{u}_{0}=\widetilde{\mathbf{f}}\text{ in }\Omega \text{,}
\label{eq2.21}
\end{equation}%
\begin{equation}
\qquad \qquad \qquad \qquad \qquad \mathbf{u}_{0}=0\text{ on }\Gamma _{1}%
\text{,}  \label{eq2.23}
\end{equation}%
\begin{equation}
\qquad \qquad \qquad \qquad \quad \mathbf{\sigma }^{0}\mathbf{n=t}\text{ on }%
\Gamma _{2}  \label{eq2.24}
\end{equation}%
where $\widetilde{\mathcal{\theta }}=\int_{\partial T}\mathcal{\theta }%
\left( y\right) d\sigma \left( y\right) $, $\widetilde{\mathbf{f}}%
=\int_{Y^{\ast }}\mathbf{f}\left( y\right) dy$ and $\mathbf{\sigma }%
^{0}=\left( \mathbf{\sigma }_{ij}^{0}\right) $.

\begin{proposition}
\label{pr3.2} The boundary value problem (\ref{eq2.21})-(\ref{eq2.24})
admits at most one weak solution in $\mathbf{V}_{0}$.
\end{proposition}

\begin{proof}
Suppose that $\mathbf{u}_{0}$ and $\mathbf{w}_{0}\in \mathbf{V}_{0}$ verify (%
\ref{eq2.21})-(\ref{eq2.24}). By the fact that the equivalent variational
problem is 
\begin{equation}
\left\{ 
\begin{array}{c}
\mathbf{u}_{0}\in \mathbf{V}_{0}\text{;}\qquad \qquad \qquad \qquad \qquad
\qquad \qquad \qquad \qquad \qquad \qquad \qquad \qquad \qquad \qquad \\ 
\sum_{i,j=1}^{N}\int_{\Omega }\mathbf{\sigma }_{ij}^{0}\left( \mathbf{u}%
_{0}\right) \mathbf{e}_{ij}\left( \mathbf{v}_{0}\right) dx+\int_{\Omega }%
\widetilde{\mathcal{\theta }}\mathbf{u}_{0}{\small \cdot }\mathbf{v}%
_{0}dx=\int_{\Omega }\widetilde{\mathbf{f}}{\small \cdot }\mathbf{v}%
_{0}dx+\int_{\Gamma _{2}}\mathbf{t{\small \cdot }v}_{0}\text{ for all }%
\mathbf{v}_{0}\in \mathbf{V}_{0}\text{,}%
\end{array}%
\right.  \label{eq2.25}
\end{equation}%
$\mathbf{z}_{0}=\mathbf{u}_{0}-\mathbf{w}_{0}$ satisfies 
\begin{equation*}
\sum_{i,j=1}^{N}\int_{\Omega }\mathbf{\sigma }_{ij}^{0}\left( \mathbf{z}%
_{0}\right) \mathbf{e}_{ij}\left( \mathbf{z}_{0}\right) dx+\int_{\Omega }%
\widetilde{\mathcal{\theta }}\mathbf{z}_{0}{\small \cdot }\mathbf{z}_{0}dx=0%
\text{.}
\end{equation*}%
Thus, using (\ref{eq1.2a}) and (\ref{eq2.21c}) we have 
\begin{equation*}
\alpha _{0}\text{mes}\left( \partial T\right) \left\Vert \mathbf{z}%
_{0}\right\Vert _{L^{2}}^{2}\leq 0
\end{equation*}%
and therefore $\mathbf{z}_{0}=0$.
\end{proof}

This leads to the following theorem.

\begin{theorem}
\label{th3.2} Suppose that the hypotheses (\ref{eq1.2a})-(\ref{eq1.2}) and (%
\ref{eq1.8a}) are satisfied. For each real $0<\varepsilon <1$, let $\mathbf{u%
}_{\varepsilon }\in \mathbf{V}_{\varepsilon }$ be defined by (\ref{eq1.3})-(%
\ref{eq1.7}) and let $\mathcal{P}^{\varepsilon }$ be the extension operator
of Proposition \ref{pr2.1}. Then, as $\varepsilon \rightarrow 0$, we have $%
\mathcal{P}^{\varepsilon }\mathbf{u}_{\varepsilon }\rightarrow \mathbf{u}%
_{0} $ in $H^{1}\left( \Omega \right) ^{N}$-weak and $\mathbf{u}_{0}$ is the
unique weak solution to (\ref{eq2.21})-(\ref{eq2.24}).
\end{theorem}

\begin{proof}
In view of the proof of Theorem \ref{th3.1}, from any given fundamental
sequence $E$ one can extract a subsequence $E^{\prime }$ such that as $%
E^{\prime }\ni \varepsilon \rightarrow 0$, we have (\ref{eq2.15a})-(\ref%
{eq2.15c}). Further (\ref{eq2.18}) holds for all $\mathbf{v=}\left( \mathbf{v%
}_{0},\mathbf{v}_{1}\right) \in \mathbb{F}_{0}^{1}$, with $\mathbf{u}=\left( 
\mathbf{u}_{0},\mathbf{u}_{1}\right) \in \mathbf{V}_{0}\times L^{2}\left(
\Omega ;\mathbf{V}_{Y}\right) $. Now, substituting (\ref{eq2.20a}) in (\ref%
{eq2.18}) and then choosing therein the $\mathbf{v}$'s such that $\mathbf{v}%
_{1}=0$, we use (\ref{eq2.21a}) and a simple computation to have (\ref%
{eq2.25}). Hence the theorem follows by Proposition \ref{pr3.2} and use of
an obvious argument.
\end{proof}

\noindent \textbf{Acknowledgements }\textit{This work is supported by the
national allocation for research modernisation of the Cameroon government.}

\end{document}